\documentclass[psamsfonts, reqno]{amsart}

\usepackage{amscd,amsmath,amssymb,amsfonts,verbatim}
\usepackage{times}
\usepackage[cmtip, all]{xy}
\usepackage{graphicx}
\usepackage[active]{srcltx}
\usepackage{color}
\usepackage{textcomp}

\definecolor{refkey}{gray}{.5}   
\definecolor{labelkey}{gray}{.5} 
\definecolor{Red}{rgb}{1,0,0}

\newcommand{\pf}{{\bf Proof : }}

\newcommand{\qedwhite}{\hfill \ensuremath{\Box}}

\newtheorem{theo}{Theorem}[section]
\newtheorem{prop}[theo]{Proposition}
\newtheorem{lem}[theo]{Lemma}
\newtheorem{cor}[theo]{Corollary}

\theoremstyle{definition}

\newtheorem{defi}[theo]{Definition}

\title{On completion of unimodular rows over polynomial extension of finitely generated rings 
over $\mathbb{Z}$ }
\author{Sampat Sharma}
\newcommand{\Addresses}{{
  \bigskip
  \footnotesize

 \textsc{Sampat Sharma, Department of Mathematics, IIT Bombay,\\ \noindent
           Main Gate Rd, IIT Area, Powai, Mumbai, Maharashtra 400076
   } \par\nopagebreak
  \textit{E-mail}: Sampat ~Sharma \texttt{<sampat@math.iitb.ac.in>; <sampat.iiserm@gmail.com>}

  \medskip

  }}

\begin{document}
\maketitle

\vskip0.15in

\subjclass 2010 Mathematics Subject Classification:{13C10, 13D15, 19B10}

\keywords {Keywords:}~ {Witt group, completion of unimodular rows}

\begin{abstract}
 In this article, we prove that if $R$ is a finitely generated ring over $\mathbb{Z}$ of dimension $d, d\geq2, \frac{1}{d!}\in R$, then any unimodular row over $R[X]$ of length $d+1$ can be mapped to a factorial row by elementary transformations.
\end{abstract}

\vskip0.50in

\begin{flushleft}
 Throughout this article we will assume $R$ to be a commutative noetherian ring with $1 \neq 0 .$
\end{flushleft}

\section{Introduction} J-P. Serre's problem on page 243 of his famous 1955 paper FAC (\cite{jps})
asked whether there were non-free projective modules over a polynomial
extension $k[X_1, \ldots, X_n]$, over a field $k$. This problem was settled
by D. Quillen (\cite{quill}) and A.A. Suslin (\cite{4}) independently in early
1976; and is now known as the Quillen--Suslin theorem. Since every f.g. projective module over $k[X_1, \ldots, X_n]$, is stably free, to determine whether projective modules are free, it is enough to determine that unimodular rows over $k[X_1, \ldots, X_n]$ are 
completable. Therefore, problem of completion of unimodular rows is a central problem in classical $K$-Theory.

\par In \cite{swantowber}, R.G. Swan and J. Towber showed that if $(a^{2}, b, c)\in Um_{3}(R)$ then it can be completed to an invertible matrix over $R.$ This result of Swan and Towber was generalised by Suslin in \cite{sus2} who showed that if 
$(a_{0}^{r!}, a_{1}, \ldots, a_{r})\in Um_{r+1}(R)$ then it can be completed to an invertible matrix. In  \cite{invent}, \cite{trans}, Ravi Rao studied the problem of completion of unimodular rows over $R[X]$, where $R$ is a local ring. Ravi Rao showed that if 
$R$ is a local ring of dimension $d, d\geq2, \frac{1}{d!}\in R$, then any unimodular row over $R[X]$ of length $d+1$ can be mapped to a factorial row by elementary transformations. 

\par In view of Suslin's local-global principle and Ravi Rao's result in $(${\cite[Theorem 2.4]{invent}}$),$ one can prove that every unimodular row of length $d+1$ over $R[X]$ where $R$ is a finitely generated ring over $\mathbb{Z}$ of dimension $d$, is completable. In this 
article, we prove that we can actually capture a factorial row in the elementary orbit of unimodular row of length $d+1.$ In particular, we prove the following :

 \begin{theo}
  Let $R$ be a finitely generated ring of dimension $d,$ $d\geq 2,$ $\frac{1}{d!}\in R$ and $v\in Um_{d+1}(R[X]).$ Then 
       $$v\overset{E_{d+1}(R[X])}{\sim} (u_{0}^{d!}, u_{1},\ldots, u_{d}) $$ for some 
            $(u_{0},\ldots, u_{d})\in Um_{d+1}(R[X]).$ 
      \end{theo} 

\par We first show that $W_{E}(R) \equiv 0$ if $R$ is a finitely generated ring over $\mathbb{Z}$ of dimension $2.$ In the next step we show that in this case, we can extract $k^{th}$ root in $W_{E}(R[X]).$ Then, using results in \cite{7}, we prove the above theorem in the 
case when $d =2.$ Next, following Roitman's degree reduction technique we can map any $v\in Um_{d+1}(R[X])$ to such a unimodular row whose last coordinate is a non-zero-divisor. Then, we conclude the result by induction hypothesis.

\section{Preliminaries}
\setcounter{theo}{0}
\par 
A row $v= (a_{0},a_{1},\ldots, a_{r})\in R^{r+1}$ is said to be unimodular if there is a $w = (b_{0},b_{1},\ldots, b_{r})
\in R^{r+1}$ with
$\langle v ,w\rangle = \Sigma_{i = 0}^{r} a_{i}b_{i} = 1$ and $Um_{r+1}(R)$ will denote
the set of unimodular rows (over $R$) of length 
$r+1$.
\par 
The group of elementary matrices is a subgroup of $GL_{r+1}(R)$, denoted by $E_{r+1}(R)$, and is generated by the matrices 
of the form $e_{ij}(\lambda) = I_{r+1} + \lambda E_{ij}$, where $\lambda \in R, ~i\neq j, ~1\leq i,j\leq r+1,~
E_{ij} \in M_{r+1}(R)$ whose $ij^{th}$ entry is $1$ and all other entries are zero. The elementary linear group 
$E_{r+1}(R)$ acts on the rows of length $r+1$ by right multiplication. Moreover, this action takes unimodular rows to 
unimodular
rows : $\frac{Um_{r+1}(R)}{E_{r+1}(R)}$ will denote set of orbits of this action; and we shall denote by $[v]$ the 
equivalence class of a row $v$ under this equivalence relation.

\par 
In $(${\cite[Theorem 3.6]{vdk1}}$),$ W. van der Kallen derives an abelian
 group structure on $\frac{Um_{d+1}(R)}{E_{d+1}(R)}$ when 
$R$ is of dimension $d,$ for all $d\geq 2.$ We will denote the group operation in this group by $\ast.$

\begin{defi} A matrix $M\in M_{r}(R)$ is said to be alternating if $M = N - N^{T}$ for some matrix 
$N\in M_{r}(R)$,  i.e. it is skew-symmetric and its diagonal entries are zero.
\end{defi}

\subsection{The elementary symplectic Witt group $W_{E}(R)$} 
If $\alpha \in M_{r}(R), \beta \in M_{s}(R)$ are matrices 
then $\alpha \perp \beta$ denotes the matrix $\begin{bmatrix}
                 \alpha & 0\\
                 0 & \beta\\
                \end{bmatrix} \in M_{r+s}(R)$. $\psi_{1}$ will denote $\begin{bmatrix}
                 0& 1\\
                -1 & 0\\
                \end{bmatrix} \in E_{2}(\mathbb{Z})$, and $\psi_{r}$ is inductively defined by $\psi_{r} = \psi_{r-1}\perp 
                \psi_{1} \in E_{2r}(\mathbb{Z})$, for $r\geq 2$.
                \par 
    If $\phi \in M_{2r}(R)$ is alternating then $\mbox{det}(\phi) = (\mbox{pf}(\phi))^{2}$ where $\mbox{pf}$ is a 
    polynomial (called the Pfaffian) in the matrix elements 
    with coefficients $\pm 1$. Note that we need to fix a sign in the choice of $\mbox{pf}$; so we insist 
    $\mbox{pf}(\psi_{r}) = 1$ 
    for all $r$. For any $\alpha \in M_{2r}(R)$ and any alternating matrix $\phi \in M_{2r}(R)$ we have 
    $\mbox{pf}(\alpha^{t}\phi
    \alpha) = \mbox{pf}(\phi)\mbox{det}(\alpha)$. For alternating matrices $\phi, \psi$ it is easy to check that 
    $\mbox{pf}(\phi \perp \psi) 
    = (\mbox{pf}(\phi))(\mbox{pf}(\psi))$.
    \par Two matrices $\alpha \in M_{2r}(R), \beta \in M_{2s}(R)$ are said to be equivalent (w.r.t. $E(R)$) if there exists  
    a matrix $\varepsilon \in SL_{2(r+s+l)}(R) \bigcap E(R)$,  such that $\alpha \perp \psi_{s+l} = \varepsilon^{t}
    (\beta \perp \psi_{r+l})\varepsilon,$ for some $l$. Denote this by $\alpha \overset{E}\sim \beta$. Then $\overset{E}\sim $ is an 
    equivalence relation; denote by $[\alpha]$ the orbit of $\alpha$ under this relation. 
    
    \par It is easy to see $(${\cite[p. 945]{7}}$)$ that $\perp$ induces the structure of an abelian group on the set of all 
    equivalence classes of alternating matrices with pfaffian $1$; this group is called elementary symplectic Witt group 
    and is denoted by $W_{E}(R)$. 

\subsection{The Vaserstein Rule} Let $R$ be a commutative ring and $v = (v_{0}, v_{1}, v_{2}), w = (w_{0}, w_{1}, w_{2}) \in Um_{3}(R)$ such that $v.w^{t} = 1.$  In \cite{7}, Vasertein associated an alternating matrix $V(v,w)$ to the pair $v,w$ : 
$$V(v,w) = \begin{bmatrix}
                 0 & v_{0} & v_{1} & v_{2}\\
                 -v_{0} & 0 & w_{2} & -w_{1}\\
                 -v_{1} & -w_{2} & 0 & w_{0}\\
                -v_{2} & w_{1} & -w_{0} & 0\\
                \end{bmatrix}\in SL_{4}(R) $$ 
with $\mbox{pf}(V(v,w)) = 1.$ In $(${\cite[Lemma 5.1]{7}}$),$ Vaserstein proved that $V(v,w)$ depands only on $v$. In $(${\cite[Lemma 5.2]{7}}$),$ the Vaserstein rule has been defined which is as follows : 
\begin{lem}$(${\cite[Lemma 5.2]{7}}$)$
\label{vasrule}  Let $R$ be a commutative ring and $v_{1} = (a_{0}, a_{1}, a_{2}), v_{2} = (b_{0}, b_{1}, b_{2})$ be two unimodular rows. Suppose $a_{0}a_{0}' + a_{1}a_{1}' + a_{2}a_{2}' = 1$ and let 
$$v_{3} = (a_{0}, (b_{1}, b_{2})\left(\begin{smallmatrix}
                 a_{1}& a_{2}\\
                -a_{2}' & a_{1}'\\
                \end{smallmatrix}\right) )\in Um_{3}(R).$$ 
Then for any $w_{1}, w_{2}, w_{3}$ such that $v_{i}.w_{i}^{t} = 1$ for $i = 1, 2, 3,$ we have 
$$[V(v_{1}, w_{1})] \perp [V(v_{2}, w_{2})] = [V(v_{3}, w_{3})]~ \mbox{in}~ W_{E}(R).$$
\end{lem}

\section{Coordinate power in Witt group} 
\setcounter{theo}{0}
Let $R$ be a finitely generated ring over $\mathbb{Z}.$ We say that an invertible alternating matrix $V$ is a coordinate $k^{th}$ power if the first row of $V$ has the form $(0, v_{1}^{k}, \ldots, v_{2r-1}).$ In this section, we will investigate under what conditions, every $[V]\in W_{E}(R[X])$ which 
is a $k^{th}$ power in $W_{E}(R[X])$ has a representative $V^{\ast}$ in its class such that $V^{\ast}$ 
is a coordinate $k^{th}$ power. We first note a result of Suslin--Vaserstein : 

\begin{theo}
\label{maintheorem} $(${\cite[Corollary 18.1, Theorem 18.2]{7}}$)$ Let $R$ be a finitely generated ring over $\mathbb{Z}$ of dimension $d, d\geq 2.$ Then $E_{d+1}(R)$ acts transitively on $Um_{d+1}(R).$ 
\end{theo}

\begin{lem}
\label{wittzero}
Let $R$ be a finitely generated ring over $\mathbb{Z}$ of dimension $2.$ Then $W_{E}(R)\equiv 0.$
\end{lem}
${\pf}$ Let $[V] \in W_{E}(R).$ In view of Theorem \ref{maintheorem}, $Um_{r}(R) = e_{1}E_{r}(R)$ for 
$r\geq 3,$ and so on applying $(${\cite[Lemma 5.3, Lemma 5.5]{7}}$)$, a few times if necessary, we can 
find an alternating matrix $W\in SL_{2}(R)$ such that $[V] = [W]$ in $W_{E}(R).$ Now, observe that $W = 
\psi_{1}.$ Thus $W_{E}(R)\equiv 0.$ 
$~~~~~~~~~~~~~~~~~~~~~~~~~~~~~~~~~~~~~~~~~~~~~~~~~~~~~~~~~~~~~~~~~~~~~~~~~~~~~~~~~~~~~~~~~~~~~~~~~~~~~~~~~
           ~~~~~\qedwhite$
           
\begin{prop}\label{root}
Let $R$ be a finitely generated ring over $\mathbb{Z}$ of dimension $2, \frac{1}{2k}\in R$ 
and $[V]\in W_{E}(R[X]).$ Then $[V] = [W]^{k}$ for some $[W]\in W_{E}(R[X]).$
\end{prop}
${\pf}$ In view of Lemma \ref{wittzero}, $W_{E}(R)\equiv 0.$ Thus we may assume that $V(0) = \psi_{r}$ 
for some $r.$ By $(${\cite[Lemma 3.1]{7}}$)$, one can find $\varepsilon \in E_{2(r+t)}(R[X])$ such that 
$$ \varepsilon^{t}(V\perp \psi_{t})\varepsilon = \psi_{r+t} + nX,$$ 
for some $t\geq 0, n\in M_{2(r+t)}(R).$ 
\par Let $\gamma = I_{2(r+t)}-\psi_{r+t}nX.$ Since $\gamma \in SL_{2(r+t)}(R[X])$, $\psi_{r+t}n$ is nilpotent, i.e. $(\psi_{r+t}n)^{l}\equiv 0$ for some $l.$ Since $\frac{1}{2k} \in R$, one can extract 
$2k^{th}$ root of $\gamma~ (= \beta^{2k})$, (see $(${\cite[Lemma 2.1]{invent}}$)$ for some $\beta \in SL_{2(r+t)}(R[X]).$ 
\par Let $\alpha = -\psi_{r+t}n.$  Thus $\gamma = I_{2(r+t)} +\alpha X,$ then 
$\beta^{k} = I_{2(r+t)} + \frac{\alpha X}{2} + \cdots.$ Since $ \varepsilon^{t}(V\perp \psi_{t})\varepsilon = \psi_{r+t} + nX $ is an alternating matrix, we get $n^{t} = -n.$  Therefore $\alpha^{t}\psi_{r+t} = \psi_{r+t}\alpha,$ whence $(\beta^{k})^{t}\psi_{r+t} = \psi_{r+t}\beta^{k}.$ Thus we have 
$$\varepsilon^{t}(V\perp \psi_{t})\varepsilon = \psi_{r+t} + nX = \psi_{r+t}\gamma = \psi_{r+t}\beta^{2k} = (\beta^{k})^{t}\psi_{r+t}\beta^{k}.$$
\par Let $W = \beta^{t}\psi_{r+t}\beta.$ Then upon applying the Whitehead lemma, one can check that 
$[V] = [W]^{k}$ in $W_{E}(R[X]).$
$~~~~~~~~~~~~~~~~~~~~~~~~~~~~~~~~~~~~~~~~~~~~~~~~~~~~~~~~~~~~~~~~~~~~~~~~~~~~~~~~~~~~~~~~~~~~~~~~~~~~~~~~~
           ~~~~~\qedwhite$
           
 \begin{prop}\label{antipodals}
 Let $R$ be a finitely generated ring over $\mathbb{Z}$ of dimension $2, \frac{1}{2}\in R$ and 
 $v = (v_{0}, v_{1}, v_{2})\in Um_{3}(R[X]).$ Then $v = (v_{0}, v_{1}, v_{2}) \underset{E_{3}(R[X])}{\sim} (-v_{0}, -v_{1}, -v_{2}).$
 \end{prop}
 ${\pf}$ Since $v = (v_{0}, v_{1}, v_{2})\in Um_{3}(R[X]),$ there exists $w = (w_{0}, w_{1}, w_{2})$ such that $v_{0}w_{0} + v_{1}w_{1} + v_{2}w_{2} = 1.$ Consider the associated alternating matrix $V = V(v,w)\in SL_{4}(R[X]).$
 \par  Since $1/2\in R,$ a famous theorem of M. Karoubi asserts that any invertible alternating matrix $V$ over the polynomial ring is stably congruent to its constant form, i.e. there exists a $\beta \in SL_{4+2l}(R[X])$, for some $l,$ such that $\beta^{T}(V\perp \psi_{l})\beta = (V(0)\perp \psi_{l}) = \psi_{l+2}.$ Last equality holds due to Lemma \ref{wittzero}.
 \par Since $\mbox{dim}(R[X]) = 3,$ by $(${\cite[Theorem 2.6]{4}}$),$ $Um_{r}(R[X]) = e_{1}E_{r}(R[X])$ for all $r\geq 5.$ Hence on applying $(${\cite[Lemma 5.3, Lemma 5.5]{7}}$)$, a few times if necessary, 
 we can find a $\beta^{\ast}\in SL_{4}(R[X])$ such that $(\beta^{\ast})^{T}V\beta^{\ast} = \psi_{2}.$  \par Let $\delta = \mbox{diagonal}(-1,1,-1,1)\in E_{4}(R[X]).$ Then $\delta^{T}\psi_{2}\delta = -\psi_{2}.$ Thus 
 \begin{equation}{\label{eq1}}\delta^{T}(\beta^{\ast})^{T}V\beta^{\ast}\delta = \delta^{T}\psi_{2}\delta = -\psi_{2} = \psi_{2}^{T} = [(\beta^{\ast})^{T}V\beta^{\ast}]^{T} = (\beta^{\ast})^{T}V^{T}\beta^{\ast}.\end{equation}
  
 Let $\sigma = (\beta^{\ast})^{T},$ then $(\sigma^{-1}\delta^{T}\sigma)V(\sigma^{-1}\delta^{T}\sigma)^{T} = -V.$ In view of $(${\cite[Corollary 1.4]{4}}$),$ $\sigma^{-1}\delta^{T}\sigma\in E_{4}(R[X]).$ Now, we get the desired upon applying $(${\cite[Theorem 10]{vas2}}$),$ to equation \ref{eq1}.
 $~~~~~~~~~~~~~~~~~~~~~~~~~~~~~~~~~~~~~~~~~~~~~~~~~~~~~~~~~~~~~~~~~~~~~~~~~~~~~~~~~~~~~~~~~~~~~~~~~~~~~~~~~
           ~~~~~\qedwhite$
 
 \begin{lem}
\label{antipodalequal} Let $R$ be a commutative ring and $v = (v_{0}, v_{1}, v_{2}) \in Um_{3}(R).$ Let us assume that $v \underset{E_{3}(R)}{\sim} (-v_{0}, -v_{1}, -v_{2}).$ Let $v^{(n)} = (v_{0}^{n}, v_{1}, v_{2})$ and let $w, w_{1} \in Um_{3}(R)$ such that 
$v.w^{t} = v^{(n)}.w_{1}^{t} = 1.$ Then 
$$[V(v, w)]^{n} = [V(v^{(n)}, w_{1})]~\mbox{in}~W_{E}(R).$$
\end{lem}
${\pf}$ We will prove it by induction on $n.$ If $n =2, $ it has been done in $(${\cite[Lemma 2.6.3]{trans}}$).$ Now assume that $n > 2.$ If $n$ is even, then write $[V(v, w)]^{n} = [V(v, w)]^{2}\perp [V(v, w)]^{n-2}.$ Now upon using induction hypothesis and 
$(${\cite[Corollary 2.6.2 (ii)]{trans}}$),$ we get the desired result.   If $n$ is odd, then write $[V(v, w)]^{n} = [V(v, w)]\perp [V(v, w)]^{n-1}.$ Now upon using induction hypothesis and 
$(${\cite[Corollary 2.6.2 (ii)]{trans}}$),$ we get the desired result.  
$~~~~~~~~~~~~~~~~~~~~~~~~~~~~~~~~~~~~~~~~~~~~~~~~~~~~~~~~~~~~~~~~~~~~~~~~~~~~~~~~~~~~~~~~~~~~~~~~~~~~~~~~~
           ~~~~~\qedwhite$

\begin{lem}
\label{coordinatepower}
 Let $R$ be a finitely generated ring over $\mathbb{Z}$ of dimension $2, \frac{1}{2}\in R$ and $v = (v_{0}, v_{1}, v_{2}), v^{(n)} = (v_{0}^{n}, v_{1}, v_{2})$ are unimodular rows of length three over $R[X].$ Let $w, w_{1}$ are such that $v.w^{t} = v^{(n)}.w_{1}^{t} = 1.$ Then 
$$[V(v, w)]^{n} = [V(v^{(n)}, w_{1})]~\mbox{in}~W_{E}(R[X]).$$
\end{lem}
${\pf}$ In view of Proposition \ref{antipodals}, $(v_{0}, v_{1}, v_{2}) \underset{E_{3}(R[X])}{\sim} (-v_{0}, -v_{1}, -v_{2}).$ Now, use Lemma \ref{antipodalequal}, to get the desired result.
$~~~~~~~~~~~~~~~~~~~~~~~~~~~~~~~~~~~~~~~~~~~~~~~~~~~~~~~~~~~~~~~~~~~~~~~~~~~~~~~~~~~~~~~~~~~~~~~~~~~~~~~~~
           ~~~~~\qedwhite$

\section{Roitman's degree reduction technique}
\setcounter{theo}{0}
Let $v\in Um_{r+1}(R[X]),$ where $r\geq \frac{d}{2} + 1,$ with $R$ a local ring of dimension $d.$ M. Roitman's argument in 
$(${\cite[Theorem 5]{roitstably}}$),$ show how one could decrease the degree of all but one (special) co-ordinate of $v.$  In this section, we follow Roitman's degree reduction technique and prove that if $R$ is finitely generated ring over $\mathbb{Z}$ of dimension 
$d$ and $v\in Um_{d+1}(R[X]),$ then  we could decrease the degree of all but one (special) co-ordinate of $v.$\\
First a few known results: 
\begin{lem}$($\cite[Corollary 2]{orbit}$)$
 \label{roit1} Let $(x_{0}, x_{1}, \ldots, x_{r})\in Um_{r+1}(R),$ $r\geq 2$ and $t$ be an element of $R$ which is invertible 
 modulo $(x_{0}, \ldots, x_{r-2})$. Then, 
 $$(x_{0}, \ldots, x_{r}) \underset{E}{\sim} (x_{0}, \ldots, x_{r-1}, t^{2}x_{r}).$$
\end{lem}

\begin{lem}$($\cite[Lemma 2]{roitstably}$)$
 \label{roit2} Let $S$ be a multiplicative subset of $R$ such that $R_{S}$ is noetherian of finite Krull dimension $d.$ Let 
 $(\overset{-}{a_{0}}, \ldots, \overset{-}{a_{r}}) \in Um_{r+1}(R_{S})$, $r>d.$ Then there exists $b_{i}\in R$ $(1\leq i\leq r)
 $ and $s\in S$ such that, $s\in R(a_{1}+b_{1}a_{0}) + \cdots + R(a_{r} + b_{r}a_{0}).$
 \end{lem}
 
 \begin{lem}$($\cite[Lemma 3]{roitstably}$)$
 \label{roit3}
  Let $f(X) \in R[X]$ have degree $n>0,$ and let $f(0)$ be a unit. Then for any $g(X) \in R[X]$ and any natural 
 number $k\geq (\mbox{degree}~g(X)-\mbox{degree}~f(X) + 1)$, there exists $h_{k}(X) \in R[X]$ of degree $< n$ such that 
 $$g(X) = X^{k}h_{k}(X)~\mbox{modulo}~(f(X)).$$
\end{lem}

\begin{lem}$($\cite[Lemma 11.1]{7}$)$
 \label{susmonic} Let $I$ be an ideal of $R[X]$ containing a monic polynomial $f$ of degree $m$. Let $g_{1}, \ldots, g_{k} \in I
 $ with degree $g_{i} < m,$ for $1\leq i\leq k.$ Assume that the coefficients of the $g_{i},$ $1\leq i\leq k,$ generate $R$. 
 Then $I$ contains a monic polynomial of degree $m-1.$
\end{lem}

\begin{lem}\label{susmoncor} Let $I$ be an ideal of $R[X]$ containing a monic polynomial. Let $m\geq 1.$ Let $g_{1}, \ldots, g_{k} \in I
 $ with degree $g_{i} < m,$ for $1\leq i\leq k.$ Assume that the coefficients of the $g_{i},$ $1\leq i\leq k,$ generate $R$. 
 Then $I$ contains a monic polynomial of every degree $\geq m-1.$
\end{lem}
${\pf}$ By Lemma \ref{susmonic}, if $I$ contains a monic polynomial of degree $n$ with $n\geq m,$ then it contains one of degree $n-1.$ Thus $I$ contains a monic polynomial of every degree $\geq m-1.$ 
$~~~~~~~~~~~~~~~~~~~~~~~~~~~~~~~~~~~~~~~~~~~~~~~~~~~~~~~~~~~~~~~~~~~~~~~~~~~~~~~~~~~~~~~~~~~~~~~~~~~~~~~~~~~~~~\qedwhite$

\begin{defi}
 Let $I$ be an ideal of a polynomial ring $R[X]$ (in one indeterminate). By $l(I)$ we denote the set consisting of $0$ and 
 all leading coefficients of $f\in I\diagdown \{0\}.$ Obviously $l(I)$ is an ideal of $R.$
\end{defi}

\begin{lem}$($H. Bass. A. Suslin$)$ $(${\cite[Chapter 3, Lemma 3.2]{lamoldbook}}$).$
\label{heightleading}
 Let $R$ be a commutative noetherian ring and $I$ be an ideal of $R[X].$ Then $\mbox{ht}_{R}l(I)\geq \mbox{ht}_{R[X]}I.$
\end{lem}

\begin{lem}
 \label{makingnzd} Let $R$ be a reduced noetherian ring and $ f = (f_{0}, f_{1}, \ldots, f_{r})\in Um_{r+1}(R[X]), r\geq 2.$ Then there 
 exists $g = (g_{0}, g_{1}, \ldots, g_{r})\in Um_{r+1}(R[X])$ in the elementary orbit of $f$ such that $l(g_{0}) = \pi,$ 
  a non-zero-divisor in $R.$
\end{lem}
${\pf}$ By $(${\cite[Corollary 9.4]{7}}$),$ there exists $\varepsilon \in E_{r+1}(R[X])$ such that 
$(f_{0}, \ldots, f_{r})\varepsilon = (g_{0}, \ldots, g_{r})$ and $\mbox{ht}(g_{1}, \ldots, g_{r}) \geq r \geq 2.$ By 
Lemma \ref{heightleading}, $\mbox{ht}_{R}l(I)\geq \mbox{ht}_{R[X]}I \geq r \geq 2,$ where $I = <g_{1}, \ldots, g_{r}>.$ In a reduced noetherian ring, the set of zero divisors is the union of minimal prime ideals,
therefore there exists $\lambda_{i} \in R$ such that $\pi = \sum_{i=1}^{r}\lambda_{i}l(g_{i}),$ a non-zero-divisor in $R.$ 
We may assume that degree $g_{0} > $ degree $g_{i}$ for $i\geq 1.$ Let degree $g_{i} = d_{i}$ for $0\leq i\leq r.$ Add 
$\lambda_{i}X^{d_{0}-d_{i}+1}g_{i}$ to $g_{0}$ to make $l(g_{0}) = \pi.$ 
$~~~~~~~~~~~~~~~~~~~~~~~~~~~~~~~~~~~~~~~~~~~~~~~~~~~~~~~~~~~~~~~~~~~~~~~~~~~~~~~~~~~~~~~~~~~~~~~~~~~~~~~~~~~~~~\qedwhite$

\begin{defi}
A polynomial $f(X)\in R[X]$ is said to be $\pi$-power monic if its highest coefficient is $\pi^{k}$ for some $k.$
\end{defi}

\begin{theo}
 \label{roittype}
 Let $R$ be a reduced finitely generated ring over $\mathbb{Z}$ of dimension $d, d\geq 3.$ Let  
 $f(X) = (f_{0}(X), f_{1}(X), \ldots, f_{d}(X))
  \in 
 Um_{d+1}(R[X]).$ Then, 
 $$(f_{0}(X), f_{1}(X), \ldots, f_{d}(X))\overset{E_{d+1}(R[X])}{\sim} (g_{0}(X), g_{1}(X), c_{2}, \ldots, c_{d})$$
 for some $(g_{0}(X), g_{1}(X), c_{2}, \ldots, c_{d})\in Um_{d+1}(R[X])$, $c_{i}\in R, 2\leq i\leq d$, and 
 $c_{d}$ is a non-zero-divisor.
\end{theo}
${\pf}$ In view of Lemma \ref{makingnzd}, we may assume that leading coefficient of $f_{0}$, $l(f_{0}) = \pi,$ a non-zero-divisor in $R.$ Let 
$\overset{-}{R} = \frac{R}{\pi R}.$ Thus, in view of  $(${\cite[Theorem 7.2]{4}}$),$   $\overset{-}{f(X)}$, can be elementarily mapped to $(1, 0, \ldots, 0).$ Upon lifting the elementary map, we  can obtain from $f(X)$, a row 
$g(X) = (g_{0}(X), g_{1}(X), \ldots, g_{d}(X))$ by elementary transformations such that 
$$g(X) \equiv (1, 0, \ldots, 0)~\mbox{mod}(R[X]\pi).$$ 
\par We can perform such transformations so that at every stage the row contains a $\pi$-power monic polynomial. Indeed, if we have to perform, e.g. the elementary transformation 
$$(h_{0},h_{1}, \ldots, h_{d}) \overset{T}{\longrightarrow} (h_{0}, h_{1}+kh_{0}, h_{2}, \ldots, h_{d})$$ 
and $h_{1}$ is $\pi$-power monic, then we replace $T$ by the following transformations
\begin{align*}
 (h_{0}, h_{1}, \ldots, h_{d})
 & \longrightarrow (h_{0}+\pi X^{n}h_{1}, h_{1}, h_{2}, \ldots, h_{d})\\
 & \longrightarrow (h_{0}+ \pi X^{n}h_{1}, h_{1}+ k(h_{0}+\pi X^{n}h_{1}) , h_{2}, \ldots, h_{d})
\end{align*}
where $n > \mbox{deg}~ h_{0}.$  We assume that 
$$(f_{0}, f_{1}, \ldots, f_{d}) \equiv (1, 0,  \ldots, 0)~\mbox{mod}(R[X]\pi)$$ 
and $f_{i}$ is $\pi$-power monic polynomial. If $i> 0$, then  replace $f_{0}$ by $f_{0} + \pi X^{n}f_{i}.$ So 
we assume that $f_{0}$ is $\pi$-power monic polynomial and $\mbox{deg}~f_{0} > 0.$ By Lemma \ref{roit3}, we assume 
$f_{i} = X^{2k}h_{i},$ where $\mbox{deg}~h_{i} < \mbox{deg}~f_{0}$ for $1\leq i\leq d.$ By Lemma \ref{roit1}, we assume 
$\mbox{deg}~f_{i} < \mbox{deg}~f_{0}$ for $1\leq i\leq d.$
\par Let $\mbox{deg}~ f_{0} = m.$ If $m = 1$, then we are done.  Assume now $m\geq 2.$ 
Let $(c_{1}, c_{2}, \ldots, c_{m(d-1)})$ be the coefficients of $1, X, \ldots, X^{m-1}$ in the polynomials 
$f_{2}(X), \ldots, f_{d}(X).$ By $($\cite[Chapter III, Lemma 1.1]{lamoldbook}$),$ the ideal generated in $R_{\pi}$ by 
$R_{\pi} \bigcap (R[X]_{\pi}\overset{-}{f_{0}} + R[X]_{\pi}\overset{-}{f_{1}})$ and the coefficients of $\overset{-}{f_{i}}$ $(
2\leq i\leq r)$ is $R_{\pi}.$ As $d\geq 3,$ $m(d-1)> \mbox{dim}R_{\pi},$ by Lemma \ref{roit2}, there exists 
$$(c_{1}', c_{2}', \ldots, c_{m(d-1)}') \equiv (c_{1}, c_{2}, \ldots, c_{m(d-1)}) ~ \mbox{mod} 
~ ((R[X]f_{0} + R[X]f_{1}) \bigcap R) $$ such that 
$R_{\pi}\overset{-}{c_{1}'} + \cdots + R_{\pi}\overset{-}{c_{m(d-1)}'} = R_{\pi}.$ Assume that we already have 
$R_{\pi}\overset{-}{c_{1}} + \cdots + R_{\pi}\overset{-}{c_{m(d-1)}} = R_{\pi}.$ By Lemma \ref{susmoncor}, the ideal 
$Rf_{0} + Rf_{2} + \cdots + Rf_{d}$ contains a $\pi$-power monic polynomial $h(X)$ of degree $m-1.$ Let leading coefficient of $h(X)$ is ${\pi}^{k}.$  Now, $f_{1}$ may be replaced with a $\pi$-power monic polynomial of degree $m-1.$ 
This does not depand on the original degree of $f_{1}.$ One first does the elementary transformations to bring the degree of $f_{1}$ down to $m-2,$ then simply adds $h$ to $f_{1}.$ In particular, we do the following elementary transformations :
\begin{align*}
 (f_{0}, f_{1}, \ldots, f_{d}) & \longrightarrow (f_{0}, {\pi}^{2k}f_{1}, \ldots, f_{d})\\
 & \longrightarrow (f_{0}, {\pi}^{2k}f_{1}- ({\pi}^{k}l(f_{1}))h,f_{2}, \ldots, f_{d})\\
& \longrightarrow (f_{0}, {\pi}^{2k}f_{1}- ({\pi}^{k}l(f_{1}))h +  h,f_{2}, \ldots, f_{d})
 \end{align*}
\par Note that in the above process, there is no change in the  the coordinate $f_{0}$.  Now, we have $\mbox{deg}~f_{0} = m$,  $\mbox{deg}~f_{1} = m-{1}.$ By Lemma \ref{roit1}, we also assume that 
$\mbox{deg}~ f_{i} < m-{1}$ for $2\leq i\leq d.$ 
\par If $m-1 = 1, $ then we are done. Assume that $m-1\geq 2.$ Let $(t_{1}, t_{2}, \ldots, t_{(m-1)(d-1)})$ be the coefficients of $1, X, \ldots, X^{m-2}$ in the polynomials 
$f_{2}(X), \ldots, f_{d}(X).$ Repeating the argument above, and using lemma \ref{susmoncor}, we can find a $\pi$-power monic polynomial $h_{1}(X)$ of degree $m-2$ in the ideal generated by $Rf_{0} + Rf_{2} + \cdots + Rf_{d}.$ Now, as above $f_{1}$ may be replaced with 
a $\pi$-power monic polynomial of degree $m-2.$ By Lemma \ref{roit1}, we also assume that 
$\mbox{deg}~ f_{i} < m-{2}$ for $2\leq i\leq d.$ Upon repeating this argument we can get $\mbox{deg}~f_{1} = 1,$ $ f_{i} \in R$ for $2\leq i\leq d.$
\par Now, it only remains to arrange that $f_{d}$ is a non-zero-divisor. Let us assume that $f_{d}$ is zero divisor, then there exists $0\neq r\in R$ such that $f_{d}.r = 0.$ Since $(f_{0}(X), f_{1}(X), f_{2}, \ldots, f_{d})\in Um_{d+1}(R[X])$, $ r \in <f_{0}(X), f_{1}(X), \ldots, f_{d-1}>.$ Therefore 
$$(f_{0}(X), f_{1}(X), f_{2}, \ldots, f_{d}) \overset{E_{d+1}(R[X])}{\sim} (f_{0}(X), f_{1}(X), f_{2}, \ldots,f_{d-1}, f_{d}+ r).$$ Since $R$ is a reduced noetherian ring, the set of zero divisors of $R$ is union of minimal prime ideals of $R$. Let $\mathfrak{P_{1}}, \ldots, \mathfrak{P_{n}}$ be finitely many minimal prime ideals of $R$. Suppose $f_{d}$ belongs to first $l$ prime ideals (this can be made possible by reenumerating them), $l<n.$ Then $r\in \cap_{i=l+1}^{n}\mathfrak{P_{i}}$ and $r$ does not belong to atleast one of the first $l$ prime ideals otherwise $r$ would be zero as $R$ is reduced ring. Therefore $f_{d} + r$ belongs to less then $l$ prime ideals. Now, inducting on number of minimal prime ideals containing $f_{d}$ we can make the last coordinate a non-zero-divisor.

$~~~~~~~~~~~~~~~~~~~~~~~~~~~~~~~~~~~~~~~~~~~~~~~~~~~~~~~~~~~~~~~~~~~~~~~~~~~~~~~~~~~~~~~~~~~~~~~~~~~~~~~~~
           ~~~~~\qedwhite$

\section{The main results} 
In this section, we prove that if $R$ is a finitely generated ring over $\mathbb{Z}$ of dimension $d, d\geq2, \frac{1}{d!}\in R$, then any unimodular row over $R[X]$ of length $d+1$ can be mapped to a factorial row by elementary transformations.
\begin{prop}
 \label{dim2prop} Let $R$ be a finitely generated ring over $\mathbb{Z}$ of dimension $2$ with $\frac{1}{2k}\in R$ and let $V\in SL_{4}(R[X])$ be an alternating
  matrix of Pfaffian $1.$ Then $[V] = [V^{\ast}]$ in $W_{E}(R[X])$ with $e_{1}V^{\ast} = (0, a ^{2k}, b, c)$, 
  and $V^{\ast} \in SL_{4}(R[X]).$ Consequently, there is a  $\gamma \in E_{4}(R[X])$ such that 
  $V = \gamma^{t}V^{\ast}\gamma.$ 
\end{prop}
${\pf}$  By Proposition \ref{root}, $[V] = [W_{1}]^{2k}$ for some $W_{1} \in W_{E}(R[X]).$ By 
$(${\cite[Theorem 2.6]{4}}$),$ $Um_{r}(R[X]) = e_{1}E_{r}(R[X])$ for $r\geq 5,$ so on applying 
$(${\cite[Lemma 5.3 and Lemma 5.5]{7}}$),$ a few times, if necessary, we can find an alternating 
matrix $W\in SL_{4}(R[X])$ 
such that $[W_{1}] = [W].$ Therefore $[V] = [W]^{2k}.$ Let $[W]^{2k} = [V^{\ast}],$ thus 
$[V] = [V^{\ast}].$
By Lemma \ref{coordinatepower}, $e_{1}V^{\ast} = (0, a^{2k}, b, c).$ 
The last statement of Proposition follows 
by applying $(${\cite[Lemma 5.3 and Lemma 5.5]{7}}$)$ and $(${\cite[Theorem 6.3]{4}}$)$
$~~~~~~~~~~~~~~~~~~~~~~~~~~~~~~~~~~~~~~~~~~~~~~~~~~~~~~~~~~~~~~~~~~~~~~~~~~~~~~~~~~~~~~~~~~~~~~~~~~~~~~~~~~~~~~\qedwhite$

\begin{theo}
\label{dim2completion}
 Let $R$ be a finitely generated ring over $\mathbb{Z}$ of Krull dimension $2$ with $\frac{1}{2k} \in R.$ Let 
 $v = (v_{0}, v_{1}, v_{2}) \in Um_{3}(R[X])$. Then there exists 
 $\varepsilon \in E_{3}(R[X])$  such that 
 $$v\varepsilon = (a^{2k}, b,c )~\mbox{for~some}~(a,b,c)\in Um_{3}(R[X]).$$
\end{theo}
${\pf}$ Choose $w = (w_{0},w_{1},w_{2})$ such that $\Sigma_{i=0}^{2}v_{i}w_{i} = 1$, and consider the alternating matrix 
$V$ with Pfaffian $1$ given by 
$$V = \begin{bmatrix}
                 0 & v_{0} & v_{1} & v_{2}\\
                 -v_{0} & 0 & w_{2} & -w_{1}\\
                 -v_{1}& -w_{2} & 0 & w_{0}\\
                 -v_{2} & w_{1} & -w_{0} & 0\\
                \end{bmatrix} \in SL_{4}(R[X]).$$ By Proposition \ref{dim2prop}, there exists an alternating 
matrix $V^{\ast}\in SL_{4}(R[X])$, with $e_{1}V^{\ast} = (0, a^{2k}, b, c),$ of Pfaffian $1$ such that 
$$[V] = [V^{\ast}].$$
Therefore there exists $\gamma \in  E_{4}(R[X])$ such that $$\gamma^{t}V\gamma = V^{\ast}.$$
By  $(${\cite[Theorem 10]{vas2}}$),$ there exists $\varepsilon \in E_{3}(R[X])$ such that 
$$v\varepsilon = (a^{2k}, b,c ).$$
$~~~~~~~~~~~~~~~~~~~~~~~~~~~~~~~~~~~~~~~~~~~~~~~~~~~~~~~~~~~~~~~~~~~~~~~~~~~~~~~~~~~~~~~~~~~~~~~~~~~~~~~~~~~~\qedwhite$ 

 \begin{theo}
      \label{sizedplus1factorial}
       Let $R$ be a finitely generated ring of dimension $d,$ $d\geq 2,$ $\frac{1}{2k}\in R$ and $v\in Um_{d+1}(R[X]).$ Then 
       $$v\overset{E_{d+1}(R[X])}{\sim} (u_{0}^{2k}, u_{1},\ldots, u_{d}) $$ for some 
            $(u_{0},\ldots, u_{d})\in Um_{d+1}(R[X]).$ 
      \end{theo}
     ${\pf}$ In view of $(${\cite[Remark 1.4.3]{invent}}$),$ we may assume that $R$ is reduced. We will prove the theorem by induction on $d.$ Case $d = 2$ is covered in Theorem \ref{dim2completion}. Let 
     $v\in Um_{d+1}(R[X])$ and $d> 2.$ Then by Theorem \ref{roittype}, there exists
       $ (g_{0}(X), g_{1}(X), c_{3}, \ldots, c_{d})\in Um_{d+1}(R[X])$
      with the property that $c_{d}\in R$ is a non-zero-divisor and 
      $c_{i}\in R, 2\leq i\leq d$ and $[ (g_{0}(X), g_{1}(X), c_{3}, \ldots, c_{d}) ] = [v].$ Let $\overset{-}
{R} = \frac{R}{(c_{d})},$  $\mbox{dim} (\overset{-}{R}) \leq d-1.$ Hence by induction hypothesis there exists 
     $\overset{-}{\varepsilon_{1}} \in E_{d}(\overset{-}{R}[X])$ such that 
     $$(\overset{-}{g_{0}(X)}, \overset{-}{g_{1}(X)}, \overset{-}{c_{3}}, \ldots, \overset{-}{c_{d-1}})
     \overset{-}{\varepsilon_{1}} = ((\overset{-}{u_{0}})^{2k}, \overset{-}{u_{1}}, \ldots, 
     \overset{-}{u_{d-1}}).$$ 
     Let $\varepsilon_{1}$ be a lift of $\overset{-}{\varepsilon_{1}}.$ Upon making appropriate elementary transformation, 
     we have 
     $$[(g_{0}(X), g_{1}(X), c_{2},\ldots, c_{d})(\varepsilon_{1}\perp 1)] = [(u_{0}^{2k}, u_{1}, \ldots, u_{d-1}, 
     c_{d})].$$ 
     Since $[ (g_{0}(X), g_{1}(X), c_{3}, \ldots, c_{d}) ] = [v],$
$$[v] = [(v_{0},\ldots,v_{d})] = [u_{0}^{2k}, u_{1}, \ldots, u_{d-1}, c_{d}].$$
$~~~~~~~~~~~~~~~~~~~~~~~~~~~~~~~~~~~~~~~~~~~~~~~~~~~~~~~~~~~~~~~~~~~~~~~~~~~~~~~~~~~~~~~~~~~~~~~~~~~~~~~~~
           ~~~~~\qedwhite$

 \begin{cor}
      \label{finalresult}
       Let $R$ be a finitely generated ring of dimension $d,$ $d\geq 2,$ $\frac{1}{d!}\in R$ and $v\in Um_{d+1}(R[X]).$ Then 
       $$v\overset{E_{d+1}(R[X])}{\sim} (u_{0}^{d!}, u_{1},\ldots, u_{d}) $$ for some 
            $(u_{0},\ldots, u_{d})\in Um_{d+1}(R[X]).$ In particular, every $v\in Um_{d+1}(R[X])$ is completable.
      \end{cor}
${\pf}$ It follws from Theorem \ref{sizedplus1factorial}, upon taking $2k =d!.$ Last statement follows from $(${\cite[Theorem 2]{sus2}}$).$
$~~~~~~~~~~~~~~~~~~~~~~~~~~~~~~~~~~~~~~~~~~~~~~~~~~~~~~~~~~~~~~~~~~~~~~~~~~~~~~~~~~~~~~~~~~~~~~~~~~~~~~~~~
           ~~~~~\qedwhite$

\begin{cor}
\label{invertingsingleelement}  Let $R$ be a finitely generated ring of dimension $d,$ $d\geq 2,$ $\frac{1}{d!}\in R,$  $s\in R$ and $v\in Um_{d+1}(R_{s}[X]).$ Then 
       $$v\overset{E_{d+1}(R_{s}[X])}{\sim} (u_{0}^{d!}, u_{1},\ldots, u_{d}) $$ for some 
            $(u_{0},\ldots, u_{d})\in Um_{d+1}(R_{s}[X]).$
\end{cor}
${\pf}$  In view of $(${\cite[Remark 1.4.3]{invent}}$),$ we may assume that $R$ is reduced. First, we observe that $R_{s} \cong \frac{R[Y]}{(sY-1)}.$ Therefore $R_{s}$ is a finitely generated ring over $\mathbb{Z}$ of dimension $\leq d.$ If $\mbox{dim}(R_{s}) = d,$ then result follows from Corollary \ref{finalresult}. If $\mbox{dim}(R_{s}) < d,$ then 
 result follows from $(${\cite[Theorem 7.2]{4}}$)$ (In this case, we can actually map it to $(1,0,\ldots,0)$ by elementary transformations.)
$~~~~~~~~~~~~~~~~~~~~~~~~~~~~~~~~~~~~~~~~~~~~~~~~~~~~~~~~~~~~~~~~~~~~~~~~~~~~~~~~~~~~~~~~~~~~~~~~~~~~~~~~~
           ~~~~~\qedwhite$

\begin{cor}
\label{invertingset}  Let $R$ be a finitely generated ring of dimension $d,$ $d\geq 2,$ $\frac{1}{d!}\in R,$  $T\subset R$ be a multiplicatively closed set  and $v\in Um_{d+1}(R_{T}[X]).$ Then 
       $$v\overset{E_{d+1}(R_{T}[X])}{\sim} (u_{0}^{d!}, u_{1},\ldots, u_{d}) $$ for some 
            $(u_{0},\ldots, u_{d})\in Um_{d+1}(R_{T}[X]).$
\end{cor}
${\pf}$  In view of $(${\cite[Remark 1.4.3]{invent}}$),$ we may assume that $R$ is reduced. Let $v = (v_{0}, v_{1}, \ldots, v_{d})\in Um_{d+1}(R_{T}[X]).$ Then there exists $s\in T$ such that $v\in Um_{d+1}(R_{s}[X]).$ In view of Corollary \ref{invertingsingleelement}, 
there exists $\varepsilon \in E_{d+1}(R_{s}[X]) \subseteq  E_{d+1}(R_{T}[X])$ such that $v\varepsilon = (u_{0}^{d!},\ldots, u_{d})$ for some 
            $(u_{0},\ldots, u_{d})\in Um_{d+1}(R_{s}[X]) \subseteq Um_{d+1}(R_{T}[X]).$
$~~~~~~~~~~~~~~~~~~~~~~~~~~~~~~~~~~~~~~~~~~~~~~~~~~~~~~~~~~~~~~~~~~~~~~~~~~~~~~~~~~~~~~~~~~~~~~~~~~~~~~~~~
           ~~~~~\qedwhite$

\medskip
\noindent
{\bf Acknowledgement:} I thank the referee for his suggestions, indicating 
a proof of Lemma \ref{susmoncor}, and for a quick review.

\Addresses

\end{document}